\documentclass[11pt,a4paper]{amsart} 
\usepackage{amssymb,amscd,amsfonts,mathrsfs, stmaryrd} 
\usepackage[mathcal]{eucal} 
\usepackage{float}
\usepackage[usenames]{color}

\theoremstyle{plain}
\newtheorem{thm}{Theorem}[section]

\newtheorem{pro}[thm]{Proposition}
\newtheorem{cor}[thm]{Corollary}

\newtheorem*{con*}{Conjecture}

\theoremstyle{remark}
\newtheorem{rem}[thm]{Remark}
\newtheorem{qun}[thm]{Question}
\newtheorem{exm}[thm]{Example}

\newtheorem{dfn}[thm]{Definition}
\newtheorem*{acknowledgements}{Acknowledgements}

\numberwithin{equation}{section}
\numberwithin{table}{section}

\newcommand{\N}{\mathbb{N}}
\newcommand{\Z}{\mathbb{Z}}
\newcommand{\Q}{\mathbb{Q}}

\newcommand{\C}{\mathbb{C}}

\newcommand{\R}{\mathbb{R}}
\newcommand{\mff}{\mathfrak{f}}

\newcommand{\gl}{\mathfrak{gl}}

\newcommand{\mfm}{\mathfrak{m}}

\newcommand{\mfp}{\mathfrak{p}}
\newcommand{\mfP}{\mathfrak{P}}

\newcommand{\lri}{\mathfrak{o}}

\newcommand{\La}{\Lambda}
\newcommand{\ol}{\overline}

\newcommand{\Gri}{\ensuremath{\mathcal{O}}}

\renewcommand{\epsilon}{\varepsilon}
\renewcommand{\phi}{\varphi}
\renewcommand{\theta}{\vartheta}
\newcommand{\mcO}{\mathcal{O}}

\DeclareMathOperator{\topo}{top}
\DeclareMathOperator{\red}{red}
\DeclareMathOperator{\Spec}{Spec}

\DeclareMathOperator{\Irr}{Irr}

\DeclareMathOperator{\rk}{rk}
\DeclareMathOperator{\SL}{SL}
\DeclareMathOperator{\GL}{GL}

\DeclareMathOperator{\SU}{SU}
\DeclareMathOperator{\GU}{GU}

\def \rarr {\rightarrow}
\def \normal {\triangleleft}
\def \T {\mathcal{T}}
\def \bfG {{\bf G}}
\def \bfH {{\bf H}}

\def \bff {{\bf f}}
\def \bfF {{\bf F}}

\def \wt {\widetilde}
\def \wh {\widehat}

\def \Fp {\ensuremath{\mathbb{F}_p}}

\def \mfo {\ensuremath{\mathfrak{o}}}
\def \mfO {\ensuremath{\mathfrak{O}}}

\def \Fq {\ensuremath{\mathbb{F}_q}}
\def \Zp  {\mathbb{Z}_p}

\begin{document}
\title[Zeta functions of groups and rings -- recent developments]{Zeta
  functions of groups and rings --\\ recent developments}

\author{Christopher Voll} \address{Fakult\"at f\"ur Mathematik, Universit\"at
  Bielefeld\\ Postfach 100131\\ D-33501 Bielefeld\\Germany}
\email{C.Voll.98@cantab.net}

\begin{abstract}
  I survey some recent developments in the theory of zeta functions
  associated to infinite groups and rings, specifically zeta functions
  enumerating subgroups and subrings of finite index or
  finite-dimensional complex representations.
\end{abstract}



\subjclass[2010]{Primary: 22E50, 20E07. Secondary: 22E55, 20F69,
  22E40, 20C15, 20G25, 11M41}

\maketitle

\thispagestyle{empty}

\section{About these notes}

Over the last few decades, zeta functions have become important tools
in various areas of asymptotic group and ring theory. With the first
papers on zeta functions of groups published barely 25 years ago, the
subject area is still comparatively young. Recent developments have
led to a wealth of results and gave rise to new perspectives on
central questions in the field. The aim of these notes is to introduce
the nonspecialist reader informally to some of these developments.

I concentrate on two types of zeta functions: firstly, zeta functions
in \emph{subgroup and subring growth} of infinite groups and rings,
enumerating finite-index subobjects. Secondly, representation zeta
functions in \emph{representation growth} of infinite groups,
enumerating finite-dimensional irreducible complex representations. I
focus on common features of these zeta functions, such as Euler
factorizations, local functional equations, and their behaviour under
base extension.

Subgroup growth of groups is a relatively mature subject area, and the
existing literature reflects this: zeta functions of groups feature in
the authoritative 2003 monograph \cite{LubotzkySegal/03} on ``Subgroup
Growth'', are the subject of the Groups St Andrews 2001
survey~\cite{duS-ennui/02} and the report~\cite{duSG/06} to the
ICM~2006. The book \cite{duSWoodward/08} contains, in particular, a
substantial list of explicit examples. Some more recent developments
are surveyed in \cite[Chapter~3]{KNV/11}.

On the other hand, few papers on representation zeta functions of
infinite groups are older than ten years. Some of the lecture notes in
\cite{KNV/11} touch on the subject. The recent
survey~\cite{Klopsch/13} on representation growth of groups
complements the current set of notes.

In this text I use, more or less as blackboxes, the theory of $p$-adic
integration and the Kirillov orbit method. The former provides a
powerful toolbox for the treatment of a number of group-theoretic
counting problems. The latter is a general method to parametrize the
irreducible complex representations of certain groups in terms of
co-adjoint orbits. Rather than to explain in detail how these tools
are employed I will refer to specific references at appropriate places
in the text. I all but ignore the rich subject of zeta functions
enumerating representations or conjugacy classes of finite groups of
Lie type; see, for instance, \cite{LiebeckShalev/05}.

\medskip

These notes grew out of a survey talk I gave at the conference Groups
St~Andrews 2013 in St~Andrews. I kept the informal flavour of the
talk, preferring instructive examples and sample theorems over the
greatest generality of the presented results. As a consequence, the
text is not a systematic treaty of the subject, but rather the result
of a subjective choice.

\section{Zeta functions in asymptotic group and ring
  theory}\label{sec:introduction}

We consider counting problems of the following general form. Let
$\Gamma$ be a -- usually infinite -- algebraic object, such as a group
or a ring, and assume that, for each $n\in\N$, we are given integers
$d_n(\Gamma)\in\N_0$, encoding some algebraic information
about~$\Gamma$. Often this data will have a profinite flavour, in the
sense that, for every $n$, there exists a finite quotient $\Gamma_n$
of $\Gamma$ such that $d_n(\Gamma)$ can be computed
from~$\Gamma_n$. In any case, we encode the sequence $(d_n(\Gamma))$
in a generating function.

\begin{dfn}
  The \emph{zeta function of $(\Gamma,(d_n(\Gamma)))$} is the
  Dirichlet generating series
  \begin{equation}\label{def:zeta}\zeta_{(d_n(\Gamma))}(s) =
    \sum_{n=1}^\infty d_n(\Gamma)n^{-s},
\end{equation}
where $s$ is a complex variable. If $(d_n(\Gamma))$ is understood from
the context, we simply write $\zeta_\Gamma(s)$
for~$\zeta_{(d_n(\Gamma))}(s)$.
\end{dfn}
In the counting problems we consider Dirichlet series often turn out
to be preferable over other generating functions, in particular if the
arithmetic function $n\mapsto d_n(\Gamma)$ satisfies some of the
following properties.

\begin{itemize}
\item[(A)] \emph{Polynomial growth}, i.e.\ the coefficients
  $d_n(\Gamma)$ -- or, equivalently, their partial sums -- have
  polynomial growth: $D_n(\Gamma):=\sum_{\nu\leq n}
  d_\nu(\Gamma)=O(n^a)$ for some $a\in\R$.
\item[(B)] \emph{Multiplicativity} in the sense of elementary number
  theory: if $n=\prod_i p_i^{e_i}$ is the prime factorization of $n$,
  then $d_n(\Gamma)=\prod_i d_{p_i^{e_i}}(\Gamma)$.
\end{itemize}
Indeed, polynomial growth implies that $\zeta_{(d_n(\Gamma))}(s)$
converges absolutely on some complex half-plane. If $d_n(\Gamma)\neq
0$ for infinitely many $n$, then the \emph{abscissa of convergence} of
$\zeta_{(d_n(\Gamma))}(s)$ is equal to
$$\alpha((d_n(\Gamma))):=\limsup_{n\rarr
  \infty}\frac{\log\sum_{\nu\leq n}D_n(\Gamma)}{\log
  n}.$$ Thus $\alpha((d_n(\Gamma)))$ gives the precise degree of
polynomial growth of the partial sums $D_n(\Gamma)$ as $n$ tends to
infinity. If the sequence $(d_n(\Gamma))$ is understood from the
context, we sometimes write $\alpha(\Gamma)$ for
$\alpha((d_n(\Gamma))$.

Multiplicativity implies that -- at least formally -- the series
\eqref{def:zeta} satisfies an \emph{Euler factorization}, indexed by
the prime numbers:
\begin{equation}\label{def:euler}
\zeta_{(d_n(\Gamma))}(s) = \prod_{p \textup{ prime}}\zeta_{(d_n(\Gamma)),p}(s),
\end{equation}
where, for a prime $p$, the
function $$\zeta_{(d_n(\Gamma)),p}(s)=\zeta_{(d_{p^i}(\Gamma))}(s)=
\sum_{i=0}^\infty d_{p^i}(\Gamma)p^{-is}$$ is called the \emph{local
  factor of $\zeta_{(d_n(\Gamma))}(s)$ at the prime $p$.} We will
later consider other Euler factorizations, indexed by places of a
number field rather than rational prime numbers, which reflect
multiplicativity features of the underlying counting problem which are
subtler than the multiplicativity of $n\mapsto d_n(\Gamma)$. In any
case, there are often \emph{rationality results} which establish that
the Euler factors are rational functions, rendering them --- at least
in principle --- amenable to computation. In practice, the study of
many (global) zeta functions of the form~\eqref{def:zeta} proceeds via
a uniform description of local factors in Euler factorizations
like~\eqref{def:euler}.

Key questions regarding zeta functions of groups and rings concern the
following:

\begin{enumerate}
\item Analytic properties regarding e.g.\ the abscissa of convergence,
  analytic continuation, natural boundaries, location and
  multiplicities of zeros and poles, residue formulae, special values
  etc.,
\item arithmetic properties of the local factors, e.g.\ rationality;
  if so, structure of numerators and denominators, special symmetries
  (functional equations) etc.,
\item the variation of these properties as $\Gamma$ varies within
  natural families of groups.
\end{enumerate}
In the sequel we survey some key results and techniques in the study
of zeta functions in the context of subgroup and subring growth
(Section~\ref{sec:subgroup}) and of representation growth
(Section~\ref{sec:representations}).

\section{Subgroup and subring growth}\label{sec:subgroup}

\subsection{Subgroup growth of finitely generated nilpotent
  groups}\label{subsec:subgroups.T}

A finitely generated group $\Gamma$ has only finitely many subgroups
of each finite index~$n$. We set, for $n\in\N$,
$$a_n(\Gamma) := \#\{H\leq \Gamma \mid |\Gamma:H| = n\}.$$ If
$s_n(\Gamma):= \sum_{\nu\leq n}a_\nu(\Gamma)=O(n^a)$ for some
$a\in\R$, then $\Gamma$ is said to be of \emph{polynomial subgroup
  growth} (\emph{PSG}). Finitely generated, residually finite groups
of PSG have been characterized as the virtually solvable groups of
finite rank; see~\cite{LubotzkyMannSegal/93}. This class of groups
includes the torsion-free, finitely generated nilpotent (or
$\T$-)groups. Let $\Gamma$ be a $\T$-group. Then the sequence
$(a_n(\Gamma))$ is multiplicative. This follows from the facts that
every finite index subgroup $H$ of $\Gamma$ contains a normal such
subgroup, and that a finite nilpotent group is isomorphic to the
direct product of its Sylow $p$-subgroups. In \cite{GSS/88},
Grunewald, Segal, and Smith pioneered the use of zeta functions in the
theory of subgroup growth of $\T$-groups. They studied the
\emph{subgroup zeta function}
$$\zeta_{\Gamma}(s):=\zeta_{(a_n(\Gamma))}(s) = \sum_{n=1}^\infty
a_n(\Gamma)n^{-s}$$ of $\Gamma$ via the Euler factorization
\begin{equation}\label{equ:euler.T}
\zeta_{\Gamma}(s) = \prod_{p \textup{ prime}}\zeta_{\Gamma,p}(s),
\end{equation}
where, for each prime $p$, the local factor at $p$ is defined
via~$\zeta_{\Gamma,p}(s)= \sum_{i=0}^\infty
a_{p^{i}}(\Gamma)p^{-is}$. One of the main result of \cite{GSS/88} is
the following fundamental theorem.

\begin{thm}\cite[Theorem~1]{GSS/88}\label{thm:rational.T}
  For all primes $p$, the function $\zeta_{\Gamma,p}(s)$ is rational
  in~$p^{-s}$, i.e.\ there exist polynomials $P_p,Q_p\in\Q[Y]$ such
  that $\zeta_{\Gamma,p}(s) = P_p(p^{-s})/Q_p(p^{-s})$. The degrees of
  $P_p$ and $Q_p$ in $Y$ are bounded.
\end{thm}
The following is by now a classical example.

\begin{exm}\cite[Proposition~8.1]{GSS/88}
  Let
\begin{equation}\label{def:heisenberg}
\bfH(\Z)=\left(\begin{matrix}1&\Z&\Z\\&1&\Z\\&&1\end{matrix}\right)
\end{equation}
 be the integral Heisenberg group. Then
\begin{equation}\label{equ:heisenberg}\zeta_{\bfH(\Z)}(s) = \zeta(s)\zeta(s-1)\zeta(2s-2)\zeta(2s-3)\zeta(3s-3)^{-1},
\end{equation}
  where $\zeta(s) = \sum_{n=1}^\infty n^{-s} = \prod_{p \textup{
      prime}}(1-p^{-s})^{-1}$ is the Riemann zeta function.
\end{exm}

It is of great interest to understand how the rational functions
giving the local zeta functions in Euler factorizations
like~\eqref{equ:euler.T} vary with the prime $p$. It is known that the
denominator polynomials $Q_p(Y)$ can be chosen to be of the form
$\prod_{i\in I}(1-p^{a_i-b_is})$, for a finite index set $I$ and
nonnegative integers $a_i$, $b_i$, all depending only
on~$\Gamma$. Computing these integers, or even just a reasonably small
set of candidates, however, remains a difficult problem. The numerator
polynomials' variation with the prime $p$ is even more mysterious. It
follows from fundamental work of du Sautoy and Grunewald that there
are finitely many varieties $V_1,\dots,V_N$ defined over $\Q$, and
rational functions $W_1(X,Y),\dots,W_N(X,Y)\in\Q(X,Y)$ such that, for
almost all primes~$p$,
\begin{equation}\label{equ:denef.T}
\zeta_{\Gamma,p}(s) = \sum_{i=1}^N | \ol{V_i}(\Fp) | W_i(p,p^{-s}),
\end{equation}
where $\ol{V_i}$ denotes the reduction of $V_i$ modulo~$p$;
cf.~\cite{duSG/00}. One may construct $\T$-groups where the numbers
$|\ol{V_i}(\Fp)|$ are not all polynomials in $p$; cf., for instance,
\cite{duS-ecII/01}. Recent results determine the degree in $Y$ of the
rational functions $P_p/Q_p\in\Q(Y)$ in Theorem~\ref{thm:rational.T}
for almost all primes~$p$; cf.\ Corollary~\ref{cor:degree}.

Variations of the sequence $(a_n(\Gamma))$ include the {normal
  subgroup sequence} $(a^\normal_n(\Gamma))$, where
$$a^\normal_n(\Gamma) := \#\{H\normal \Gamma \mid |\Gamma:H|=n\}.$$
It gives rise to the \emph{normal (subgroup) zeta function}
$$\zeta^{\triangleleft}_{\Gamma}(s) :=
\zeta_{(a^\triangleleft_n(\Gamma))}(s) = \sum_{n=1}^\infty
a^\normal_n(\Gamma)n^{-s}$$ of $\Gamma$. It also has an Euler
factorization whose factors are rational in~$p^{-s}$ and, in
principle, given by formulae akin to~\eqref{equ:denef.T}.  The {normal
  zeta function} of the integral Heisenberg group
(cf.~\eqref{def:heisenberg}), for example, is equal to
$$\zeta^{\triangleleft}_{\bfH(\Z)}(s) =\zeta(s)\zeta(s-1)\zeta(3s-2) =
\prod_{p \textup{
    prime}}\frac{1}{(1-p^{-s})(1-p^{1-s})(1-p^{2-3s})};$$
cf.~\cite[Section 8]{GSS/88}.

It is interesting to ask how subgroup zeta functions of $\T$-groups,
or their variations, vary under base extension. Given a number field
$K$ with ring of integers $\mcO$ one may consider, for instance, the
$\T$-group $\bfH(\mcO)$ of upper-unitriangular $3\times3$-matrices
over~$\mcO$. Then
\begin{equation}\label{equ:euler.heisenberg}
  \zeta^{\triangleleft}_{\bfH(\mcO)}(s) = \prod_{p \textup{
      prime}}\zeta^{\triangleleft}_{\bfH(\mcO),p}(s).
\end{equation}
The following result extends parts of~\cite[Theorem~2]{GSS/88} and
makes it more precise.

\begin{thm}\cite{SV1/13}\label{thm:BSV}
  For every $r\in\N$ and every finite family $\bff =
  (f_1,\dots,f_r)\in\N^r$, there exist \emph{explicitly given}
  rational functions $W_\bff(X,Y)\in\Q(X,Y)$, such that the following
  hold.
\begin{enumerate}
\item If $p$ is a prime which is unramified in $K$ and decomposes in
  $K$ as $p\mcO=\prod_{i=1}^r \mfP_i$ for prime ideals $\mfP_i$ of
  $\Gri$ with inertia degrees $f_i = \log_p|\mcO:\mfP_i|$ for
  $i=1,\dots,r$, then
  \begin{equation*} \zeta^{\triangleleft}_{\bfH(\mcO),p}(s) =
    W_\bff(p,p^{-s}).
  \end{equation*}
\item Setting $d = | K:\Q| = \sum_{i=1}^r f_i$, we have
\begin{equation}\label{equ:funeq.sv}
 W_\bff(X^{-1},Y^{-1}) = (-1)^n X^{\binom{3d}{2}}Y^{5d}W_\bff(X,Y).
\end{equation}
\end{enumerate}
\end{thm}

The proof of Theorem~\ref{thm:BSV} is essentially combinatorial. In
the case that $p$ splits completely, i.e.\ $\bff=(1,\dots,1)$, it
proceeds by organizing the infinite sums defining the local zeta
functions as sums indexed by pairs of partitions $(\lambda,\mu)$, each
of at most $n$ parts, where $\lambda$ dominates~$\mu$. We further
partition the infinite set of such pairs into $C_n =
\frac{1}{n+1}\binom{2n}{n}$ (the $n$-th Catalan number) parts, indexed
by the Dyck words of length~$2n$, determined by the ``overlap''
between $\lambda$ and~$\mu$. This subdivision by Dyck words is
suggested by a simple lemma, attributed to Birkhoff, that determines
the numbers of subgroups of type $\mu$ in a finite abelian $p$-group
of type $\lambda$. For each fixed Dyck word, we express the
corresponding partial sum of the local zeta function in terms of
natural generalizations of combinatorially defined generating
functions, first studied by Igusa (cf.~\cite[Theorem~4]{Voll/05}) and
Stanley \cite{Stanley/91}. Remarkably, a functional equation of the
form~\eqref{equ:funeq.sv} is already satisfied by each of the $C_n$
partial sums.  If $p$ does not split completely, the strategy above
still works after some moderate modification.

The functional equation~\eqref{equ:funeq.sv} reflects the Gorenstein
property of certain face rings. That such a functional equation holds
for \emph{almost all} primes $p$ follows
from~\cite[Theorem~B]{Voll/10}; that it holds in fact for all
unramified primes is additional information. Note that $3d =
h(\bfH(\mcO))$ and $5d = h(\bfH(\mcO)) + h(\bfH(\mcO)/Z(\bfH(\mcO)))$,
the sums of the Hirsch lengths of the nontrivial quotients by the
terms of the upper central series of~$\bfH(\mcO)$. Here, given a
$\T$-group $G$, we write $h(G)$ for the Hirsch length of $G$,
i.e.\ the number of infinite cyclic factors in a decomposition series
of~$G$.

Formulae for the Euler factors in \eqref{equ:euler.heisenberg} indexed
by primes which are nonsplit (but possibly ramified) in $K$ are given
in~\cite{SV2/14}.

\subsection{Subring growth of additively finitely generated rings}

By a \emph{ring} we shall always mean a finitely generated,
torsion-free abelian group, together with a bi-additive multiplication
-- not necessarily associative, commutative, or unital. Examples of
such rings include $\Z^d$ (e.g.\ with null-multiplication or with
componentwise multiplication), the rings of integers in number fields,
and {Lie rings}, that is rings with a multiplication (or ``Lie
bracket'') which is alternating and satisfies the Jacobi
identity. Examples of Lie rings include ``semi-simple'' matrix rings
such as $\mathfrak{sl}_N(\Z)$ and the Heisenberg Lie ring
$$\mathfrak{h}(\Z)
=\left(\begin{matrix}0&\Z&\Z\\&0&\Z\\&&0\end{matrix}\right),$$ with
Lie bracket induced from $\gl_3(\Z)$.

The subring sequence of a ring $\Lambda$ is $(a_n(\Lambda))$,
where $$a_n(\Lambda):= \#\{H\leq \Gamma \mid |\Gamma:H| = n\}.$$ It is
encoded in the \emph{subring zeta function} of $\Lambda$, that is the
Dirichlet generating series
$$\zeta_{\Lambda}(s) = \zeta_{(a_n(\Lambda))}(s) = \sum_{n=1}^\infty
a_n(\Lambda)n^{-s}.$$ In contrast to the case of subgroup growth,
polynomial growth requires no assumption on the multiplicative
structure: indeed, the null-multiplication on $\Z^d$ yields a trivial
polynomial upper bound on $s_n(\Lambda):=\sum_{\nu\leq
  n}a_\nu(\Lambda)$.  Also, multiplicativity of the subring growth
function $n\mapsto a_n(\Lambda)$ follows from the Chinese Reminder
Theorem. Consequently, the subring zeta function of $\Lambda$
satisfies the following Euler factorization:
$$\zeta_{\Lambda}(s) = \prod_{p \textup{ prime}}\zeta_{\Lambda,p}(s).$$
Many of the structural results for local zeta functions of $\T$-groups
have analogues in the setting of zeta functions of rings. One example
is the following.
\begin{thm}\cite[Theorem~3.5]{GSS/88}\label{thm:rational}
  For all primes $p$, the function $\zeta_{\Lambda,p}(s)$ is rational
  in~$p^{-s}$, i.e.\ there exist polynomials $P_p,Q_p\in\Q[Y]$ such
  that $\zeta_{\Lambda,p}(s) = P_p(p^{-s})/Q_p(p^{-s})$. The degrees
  of $P_p$ and $Q_p$ in $Y$ are bounded.
\end{thm}

As in the context of subgroup growth of $\T$-groups, one also
considers variations such as {ideal growth} of rings. The \emph{ideal
  zeta function} of a ring enumerates its ideals of finite additive
index. These zeta functions, too, enjoy Euler factorizations indexed
by the rational primes. A rationality result analogous to
Theorem~\ref{thm:rational} holds for the local factors. From this
perspective one recovers, for example, the classical Dedekind zeta
function of a number field, enumerating ideals of finite index in the
number field's ring of integers.

In fact, the study of subgroup zeta functions of $\T$-groups as
outlined in Section~\ref{subsec:subgroups.T} may -- to a large extent
-- be reduced to the study of subring zeta functions of nilpotent Lie
rings. Indeed, a key tool in the analysis of \cite{GSS/88} is a
linearization technique: the Mal'cev correspondence associates to each
$\T$-group $\Gamma$ a nilpotent Lie ring $\Lambda(\Gamma)$, that is a
Lie ring whose additive group is isomorphic to $\Z^d$, where
$d=h(\Gamma)$ is the Hirsch length of~$\Gamma$, which is nilpotent
with respect to the Lie bracket; see \cite[Section~4]{GSS/88} for
details on the Mal'cev correspondence and its consequences for zeta
functions of $\T$-groups. One of these consequences is the fact that,
for almost all primes $p$,
\begin{equation}\label{equ:malcev}, 
  \zeta_{\Gamma,p}(s) = \zeta_{\Lambda(\Gamma),p}(s);
\end{equation}
cf.~\cite[Theorem~4.1]{GSS/88}. In nilpotency class at most $2$ this
equality holds for all primes~$p$. The formula \eqref{equ:heisenberg},
for instance, coincides with the subring zeta function of the
Heisenberg Lie ring~$\mathfrak{h}(\Z)=\Lambda(\bfH(\Z))$.

Maybe it is due to connections to subgroup growth like the ones just
sketched that the study of subring growth has long focussed on
\emph{Lie} rings. The following example does not arise in this
context.

\begin{exm}
  Let $\mcO$ be the ring of integers in a number field $K$ and, for
  $n\in\N$, let $b_n(\mcO)$ denote the number of subrings of $\mcO$ of
  index $n$, containing~$1\in\mcO$. The resulting zeta function
  $\zeta_{(b_n(\mcO))}(s)$ may be called the \emph{order zeta
    function} $\eta_K(s)$ of~$K$. The function $\eta_K$ has an Euler
  factorization indexed by the rational primes, though -- in contrast
  to the Dedekind zeta function $\zeta_K(s)$ -- not generally by the
  prime ideals of~$\mcO$. Clearly $\eta_\Q=1$. If $d=|K:\Q| = 2$, then
  $\eta_K(s) = \zeta(s)$, the Riemann zeta function. For $d=3$ it is
  known that
$$\eta_K(s) = \frac{\zeta_K(s)}{\zeta_K(2s)}\zeta(2s)\zeta(3s-1);$$
  see~\cite{DatskovskyWright/88}. For $d=4$, Nakagawa computes in
  \cite{Nakagawa/96} the Euler factors $\eta_{K,p}(s)$, where $p$
  ranges over the primes with arbitrary but fixed decomposition
  behaviour in~$K$.  Remarkably, the resulting formulae are rational
  functions in $p$ and $p^{-s}$ though not, in general, expressible in
  terms of translates of local Dedekind zeta functions. It is
  interesting to establish whether this uniformity on sets of primes
  with equal decomposition behaviour is a general phenomenon. Of
  particular interest is the case of primes with split totally in~$K$,
  i.e.\ primes $p$ such that $p\mcO = \mfP_1\cdots\mfP_d$, where
  $\mfP_1,\dots,\mfP_d$ are pairwise distinct prime ideals of $\mcO$
  with trivial residue field extension. For such primes, it is not
  hard to see that
$$\eta_{K,p}(s) = \zeta_{\Z^{d-1},p}(s),$$
where we consider $\Z^{d-1}$ as a ring with componentwise
multiplication.

\begin{thm}\cite{Nakagawa/96} and
  \cite[Proposition~6.3]{Liu/07}. Consider $\Z^3$ as ring with
  componentwise multiplication. Then $\zeta_{\Z^3}(s) = \prod_{p
    \textup{ prime}}\zeta_{\Zp^3,p}(s)$, where, setting $t=p^{-s}$, we
  have
\begin{multline}
\zeta_{\Z^3,p}(s) = \label{equ:liu}\\ \frac{1 + 4t + 2t^2 + (4p -
  3)t^3 + (5p - 1)t^4 + (p^2-5p)t^5 + (3p^2-4p)t^6 - 2p^2 t^7 -
  4p^2t^8 - p^2 t^9}{(1-t)^2(1-p^2t^4)(1-p^3t^6)}.
\end{multline}
\end{thm}
\end{exm}

The evidence available for $d\leq 3$ suggests a positive answer to the
following question.
\begin{qun}
  Do there exist rational functions $W_d(X,Y)\in\Q(X,Y)$, for
  $d\in\N$, such that, for all primes $p$,
$$\zeta_{\Zp^d,p}(s) = W_d(p,p^{-s})?$$
\end{qun}

Local zeta functions such as the ones given in~\eqref{equ:liu} exhibit
a curious palindromic symmetry under inversion of~$p$. This is no
coincidence, as the following result shows.

\begin{thm}\cite[Theorem~A]{Voll/10}\label{thm:funeq}
  Let $\Lambda$ be a ring with $(\Lambda,+)\cong \Z^d$. Then, for
  almost all primes $p$,
\begin{equation}\label{equ:funeq}
\zeta_{\Lambda,p}(s)|_{p \rarr p^{-1}} = (-1)^d p^{\binom{d}{2}-ds}\zeta_{\Lambda,p}(s).
\end{equation}
\end{thm}

\begin{cor}\label{cor:degree}
For almost all primes $p$, $\deg_{p^{-s}}(\zeta_{\Lambda,p}(s)) = -d$.
\end{cor}
Via the Mal'cev correspondence, Theorem~\ref{thm:funeq} yields an
analogous statement for almost all of the local factors
$\zeta_{\Gamma,p}(s)$ of the subgroup zeta function $\zeta_\Gamma(s)$
of a $\T$-group~$\Gamma$; cf.~\eqref{equ:malcev}. There are analogous
results giving functional equations akin to~\eqref{equ:funeq} for
ideal zeta functions of $\T$-groups -- or equivalently, again by the
Mal'cev correspondence, nilpotent Lie rings of finite additive rank --
of nilpotency class at most~$2$. There are, however, examples of
$\T$-groups of nilpotency class $3$ whose local normal subgroup zeta
functions do not satisfy functional equations like~\eqref{equ:funeq};
cf.~\cite[Theorem~1.1]{duSWoodward/08}.

Other variants of subgroup zeta functions of $\T$-groups which have
been studied include those encoding the numbers of finite-index
subgroups whose profinite completion is isomorphic to the one of the
ambient group. These \emph{pro-isomorphic zeta functions} also enjoy
Euler product decompositions, indexed by the rational primes, whose
factors are rational functions. It is an interesting open problem to
characterise the $\T$-groups for which these local factors satisfy
functional equations comparable to~\eqref{equ:funeq}. For positive
results in this direction see \cite{duSLubotzky/96, Berman/11}. An
example of a $\T$-group (of nilpotency class $4$ and Hirsch length
$25$) whose pro-isomorphic zeta function's local factors do not
satisfy such functional equations was recently given in
\cite{BermanKlopsch/14}.

\subsection{Taking the limit $p\rightarrow 1$: reduced and topological
  zeta functions of groups and rings}
Numerous mathematical concepts, theorems, and identities allow natural
$q$-{analogues}. Featuring an additional parameter $q$ -- often
interpreted as a prime power --, these analogues return the original
object upon setting $q=1$. Examples include the Gaussian $q$-binomial
coefficients, generalizing classical binomial coefficients and Heine's
basic hypergeometric series, generalizing ordinary hypergeometric
series.

An idea that only recently took hold in the theory of zeta functions
of groups and rings is to interpret local such zeta functions as
``$p$-analogues'' of certain limit objects as $p\rarr 1$ and to
investigate the limit objects with tools from combinatorics or
commutative algebra.

\subsubsection{Reduced zeta functions} One way to make this idea
rigorous leads, for instance, to the concept of the \emph{reduced zeta
  function} $\zeta_{\Lambda,\red}(t)$ of a ring $\Lambda$. Informally,
this rational function in a variable $t$ over the rationals is
obtained by setting $p=1$ in the coefficients of the $p$-adic subring
zeta function of $\Lambda$, considered as a series in~$t=p^{-s}$;
formally, it arises by specializing the coefficients of the motivic
zeta function associated to $\Lambda$ via the Euler characteristic;
cf.~\cite{duSLoeser/04, Evseev/09}. Under some very restrictive
conditions on $\Lambda$, the reduced zeta function
$\zeta_{\Lambda,\red}(t)$ is known to enumerate the integral points of
a rational polyhedral cone. In the language of commutative algebra
this means that $\zeta_{\Lambda,\red}(t)$ is the Hilbert series of an
affine monoid algebra attached to a Diophantine system of linear
inequalities. For general rings a somewhat more multifarious picture
seems to emerge, as the following example indicates.

\begin{exm}\label{exm:Z3.redu}
  Consider $\Lambda=\Z^3$ as a ring with componentwise
  multiplication. Heuristically, setting $p=1$ in~\eqref{equ:liu} we obtain 
$$\zeta_{\Z^3,\red}(t) = \frac{1 + 5t + 6t^2 + 3t^3 + 6t^4 + 5t^5 +
  t^6}{(1-t)(1-t^2)(1-t^6)}.$$ Intriguingly, this rational function is
{not} the generating function of a polyhedral cone, but does exhibit
some tell-tale signs of the Hilbert series of a graded Cohen-Macaulay
(even Gorenstein) algebra of dimension~$3$.
\end{exm}

\subsubsection{Topological zeta functions}
Topological zeta functions offer another way to define a limit as
$p\rarr 1$ of families of $p$-adic zeta functions. They were first
introduced in the realm of Igusa's $p$-adic zeta function as
singularity invariants of hypersurfaces
\cite{DenefLoeser/92}. Informally, the topological zeta function is
the leading term of the expansion of the $p$-adic zeta function in
$p-1$. Formally, it may be obtained by specialising the motivic zeta
function; cf.~\cite{DenefLoeser/98}. Whereas the latter lives in the
power series ring over a certain completion of a localization of a
Grothendieck ring of algebraic varieties, the topological zeta
function is just a rational function in one variable $s$, say, over
the rationals. The topological zeta function
$\zeta_{\Lambda,\topo}(s)$ of a ring was introduced
in~\cite{duSLoeser/04}.

\begin{exm}
  The topological zeta function of $\Z^3$ (cf.\
  Example~\ref{exm:Z3.redu}) is $$\zeta_{\Z^3,\topo}(s) =
\frac{9s-1}{s^2(2s-1)^2}.$$
\end{exm}
In \cite{RossmannI/14} Rossmann develops an effective method for
computing topological zeta functions associated to groups, rings, and
modules. It is built upon explicit convex-geometric formulae for a
class of $p$-adic integrals under suitable non-degeneracy conditions
with respect to associated Newton polytopes. This method yields
examples of explicit formulae for topological zeta functions of
objects whose $p$-adic zeta functions are well out of computational
reach. For a number of intriguing conjectures about arithmetic
properties of topological zeta functions
see~\cite[Section~8]{RossmannI/14}. Rossmann implemented his algorithm
in \cite{RossmannZeta/14} and explained it in detail
in~\cite{RossmannII/14}.

\section{Representation growth}\label{sec:representations}

Let $\Gamma$ be a group. Consider, for $n\in\N$, the set
$\Irr_n(\Gamma)$ of $n$-dimensional irreducible complex
representations of $\Gamma$ up to isomorphism. If $\Gamma$ has
additional structure, we restrict our attention to representations
respecting this structure. For instance, if $\Gamma$ is a topological
group, we only consider continuous representations. The group $\Gamma$
is called (\emph{representation}) \emph{rigid} if
$r_n(\Gamma):=\#\Irr_n(\Gamma)$ is finite for all $n$. In this case,
the Dirichlet generating series
$$\zeta^{\textup{irr}}_\Gamma(s):=\zeta_{(r_n(\Gamma))}(s) = \sum_{n=1}^\infty r_n(\Gamma)n^{-s}$$ is called the
\emph{representation zeta function} of $\Gamma$.  We discuss several
classes of groups whose representation zeta functions (or natural
variants thereof) have recently attracted attention. These are
\begin{enumerate}
\item finitely generated nilpotent groups,
\item arithmetic groups in characteristic $0$,
\item algebraic groups,
\item compact $p$-adic analytic groups,
\item iterated wreath products and branch groups.
\end{enumerate}
Throughout, let $K$ be a number field with ring of integers
$\mcO=\mcO_K$. We write
$$\zeta_K(s)=\sum_{I\triangleleft
  \mcO}|\mcO:I|^{-s}=\prod_{\mfp\in\Spec{\mcO}}(1-|\mcO/\mfp|^{-s})^{-1}$$
for the Dedekind zeta function of $K$. Note that
$\zeta_\Q(s)=\zeta(s)$.  By representations we will always mean
complex representations.

\subsection{Finitely generated nilpotent groups}

Let $\Gamma$ be a $\T$-group. Unless $\Gamma$ is trivial, the sets
$\Irr_n(\Gamma)$ are not all finite. Indeed, a nontrivial $\T$-group
surjects onto the infinite cyclic group and thus has infinitely many
one-dimensional representations. We therefore consider
finite-dimensional representations up to twists by one-dimensional
representations. More precisely, two representations
$\rho_1,\rho_2\in\Irr_n(\Gamma)$ are said to be
\emph{twist-equivalent} if there exists $\chi\in\Irr_1(\Gamma)$ such
that $\rho_1$ is isomorphic to $\rho_2\otimes\chi$. The numbers
$\wt{r}_n(\Gamma)$ of isomorphism classes of irreducible, complex
$n$-dimensional representations of $\Gamma$ are all finite;
cf.~\cite[Theorem~6.6]{LubotzkyMagid/85}. We define the
\emph{representation zeta function} of $\Gamma$ to be the Dirichlet
generating series
$$\zeta^{\wt{\textup{irr}}}_\Gamma(s) := \zeta_{(\wt{r}_n(\Gamma))}(s) = \sum_{n=1}^\infty\wt{r}_n(\Gamma)n^{-s}.$$
The coefficients $\wt{r}_n(\Gamma)$ grow polynomially, so
$\zeta^{\wt{\textup{irr}}}_\Gamma(s)$ converges on some complex
right-half plane. The precise abscissa of convergence of
$\zeta^{\wt{\textup{irr}}}_\Gamma(s)$ is an interesting invariant
of~$\Gamma$.

The function $n \mapsto \wt{r}_n(\Gamma)$ is multiplicative, which
yields the Euler factorization
\begin{equation}\label{equ:euler.ztwirr}
  \zeta^{\wt{\textup{irr}}}_\Gamma(s) = \prod_{p \textup{ prime}}\zeta^{\wt{\textup{irr}}}_{\Gamma,p}(s),
\end{equation}
where, for a prime~$p$, the local factor
$\zeta^{\wt{\textup{irr}}}_{\Gamma,p}(s) =
\sum_{i=0}^\infty\wt{r}_{p^i}(\Gamma)(p^{-s})^i$ enumerates
twist-isoclasses of representations of $\Gamma$ of $p$-power
dimension.

\begin{thm}\cite[Theorem~8.5]{HrushovskiMartinRideauCluckers/14}
  For all primes $p$, the function
  $\zeta^{\wt{\textup{irr}}}_{\Gamma,p}(s)$ is rational in~$p^{-s}$,
  i.e.\ there exist polynomials $P_p,Q_p\in\Q[Y]$ such that
  $\zeta^{\wt{\textup{irr}}}_{\Gamma,p}(s) =
  P_p(p^{-s})/Q_p(p^{-s})$. The degrees of $P_p$ and $Q_p$ in $Y$ are
  bounded.
\end{thm}
The proof uses model-theoretic results on definable equivalence
classes. We illustrate this important rationality result with a simple
but instructive example.

\begin{exm}\label{exa:heisenberg.ztwirr}
  Consider the integral Heisenberg group $\bfH(\Z)$;
  cf.~\eqref{def:heisenberg}. Then
\begin{equation}\label{equ:euler.ztwirr.Heisenberg.Q}
\zeta^{\wt{\textup{irr}}}_{\bfH(\Z)}(s) = \sum_{n=1}^\infty \phi(n)n^{-s} =
\frac{\zeta(s-1)}{\zeta(s)}=\prod_{p \text{ prime}}\frac{1-p^{-s}}{1-p^{1-s}},
\end{equation}
where $\phi$ is the Euler totient function; cf.~\cite{NunleyMagid/89}.

It turns out that the formula in \eqref{equ:euler.ztwirr.Heisenberg.Q}
behaves uniformly under some base extensions, as we shall now explain.
Consider, for example, the Heisenberg group $\bfH(\mcO)$ over $\mcO$,
i.e.\ the group of upper-unitriangular $3\times 3$-matrices
over~$\mcO$. Then
\begin{equation}\label{equ:euler.ztwirr.Heisenberg.no.field}
\zeta^{\wt{\textup{irr}}}_{\bfH(\mcO)}(s) = \frac{\zeta_{K}(s-1)}{\zeta_K(s)}=\prod_{\mfp\in\Spec{\mcO}}\frac{1-|\mcO/\mfp|^{-s}}{1-|\mcO/\mfp|^{1-s}}.
\end{equation} 
For quadratic number fields this was proven by Ezzat
in~\cite{Ezzat/14}. The general case follows
from~\cite[Theorem~B]{StasinskiVoll/14}.

Each factor
$\zeta^{\wt{\textup{irr}}}_{\bfH(\mcO),\mfp}(s):=\frac{1-|\mcO/\mfp|^{-s}}{1-|\mcO/\mfp|^{1-s}}$
of the Euler factorization
\eqref{equ:euler.ztwirr.Heisenberg.no.field} is interpretable as a
representation zeta function associated to a pro-$p$ group. Indeed,
for $\mfp\in\Spec \mcO$, we denote by $\mcO_\mfp$ the completion of
$\mcO$ at~$\mfp$. Then
$\zeta^{\wt{\textup{irr}}}_{\bfH(\mcO),\mfp}(s)$ is equal to the zeta
function $\zeta^{\wt{\textup{irr}}}_{\bfH(\mcO_\mfp)}(s)$ of the
pro-$p$ group
$$\bfH(\mcO_\mfp)=\left(\begin{matrix}1&\mcO_\mfp&\mcO_\mfp\\&1&\mcO_\mfp\\&&1\end{matrix}\right),$$
enumerating \emph{continuous} irreducible representations of
$\bfH(\mcO_\mfp)$ up to twists by \emph{continuous} one-dimensional
representations. We note the following features
of~$\zeta^{\wt{\textup{irr}}}_{\bfH(\mcO)}(s)$.
\begin{enumerate}
\item Whilst the Euler factorization
  \eqref{equ:euler.ztwirr.Heisenberg.Q} illustrates the general
  factorization~\eqref{equ:euler.ztwirr}, the factorization
  \eqref{equ:euler.ztwirr.Heisenberg.no.field} is finer than
  \eqref{equ:euler.ztwirr}. In fact, for each rational
  prime~$p$, $$\zeta^{\wt{\textup{irr}}}_{\bfH(\mcO),p}(s) =
  \prod_{\mfp | p\mcO}
  \zeta^{\wt{\textup{irr}}}_{\bfH(\mcO_\mfp)}(s).$$
\item The factors of the ``fine'' Euler factorization
  \eqref{equ:euler.ztwirr.Heisenberg.no.field} are indexed by the
  nonzero prime ideals $\mfp$ of $\mcO$, and are each given by a
  rational functions in $q^{-s}$, where $q=|\mcO/\mfp|$ denotes the
  residue field cardinality.
\item Each factor of the Euler factorization
  \eqref{equ:euler.ztwirr.Heisenberg.no.field} satisfies the
  functional equation
$$\left.\zeta^{\wt{\textup{irr}}}_{\bfH(\mcO),\mfp}(s)\right|_{q\rarr q^{-1}}=
\left.\frac{1-q^{-s}}{1-q^{1-s}}\right|_{q\rarr q^{-1}} =
\frac{1-q^s}{1-q^{-1+s}}=q\,\zeta^{\wt{\textup{irr}}}_{\bfH(\mcO),\mfp}(s).$$
\end{enumerate}
\end{exm}
As we shall see, all of these points are special cases of general
phenomena. 

We consider in the sequel families of groups obtained from Lie
lattices. Let, more precisely, $\Lambda$ be an $\mcO$-Lie lattice,
i.e.\ a free and finitely generated $\mcO$-module, together with an
antisymmetric, bi-additive
form~$[\,,]:\Lambda\times\Lambda\rarr\Lambda$, called `Lie bracket',
which satisfies the Jacobi identity. Assume further that $\Lambda$ is
nilpotent with respect to $[\,,]$ of class~$c$, and let $\Lambda'$
denote the derived Lie lattice $[\Lambda,\Lambda]$. If $\Lambda$
satisfies $\Lambda' \subseteq c!\Lambda$, then it gives rise to a
unipotent group scheme $\bfG_\Lambda$ over $\mcO$, via the Hausdorff
series as we shall now explain. The Hausdorff series $F(X,Y)$ is a
formal power series in two noncommuting variables $X$ and $Y$, with
rational coefficients. The Hausdorff formula gives an expression for
this series in terms of Lie terms:
\begin{equation}\label{equ:hausdorff}
F(X,Y) = X + Y + \frac{1}{2}[X,Y] + \frac{1}{12}\left( [[X,Y],Y] - [[X,Y],X]\right) + \dots,
\end{equation}
where $[A,B]:=AB-BA$. See, e.g., \cite[Chapter~I, Section~7.4]{KNV/11}
for further details on the Hausdorff series.

For an $\mcO$-algebra $R$, let $\Lambda(R) := \Lambda\otimes_\mcO
R$. The assumption that $\Lambda' \subseteq c!\Lambda$ allows one to
define on the set $\Lambda(R)$ a group structure $*$ by setting, for
$x,y\in\Lambda(R)$,
$$x*y := F(x,y), \quad x^{-1} = -x.$$ Note that, by the nilpotency of
$\Lambda$, the Hausdorff formula \eqref{equ:hausdorff} yields an
expression for $x*y$ as a linear combination of Lie terms in $x$
and~$y$.  In this way one obtains a unipotent group scheme
$\bfG_\Lambda$ over $\mcO$, representing the functor $R \mapsto
(\Lambda(R),*)$. In nilpotency class $c=2$, one may define the group
scheme $\bfG_\Lambda$ directly and avoiding the condition
$\Lambda'\subseteq c!\Lambda$;
cf.~\cite[Section~2.4]{StasinskiVoll/14}.

By taking rational points of $\bfG_\Lambda$ we obtain a multitude of
groups, all originating from the same global Lie
lattice~$\Lambda$. The group $\bfG_\Lambda(\mcO)$ of $\mcO$-rational
points, for instance, is a $\T$-group of Hirsch length
$\rk_\Z(\mcO)\rk_\mcO(\Lambda)$. By considering the
$\Gri_\mfp$-rational points of $\bfG_\Lambda$ for a nonzero prime
ideal $\mfp$ of $\Gri$, we obtain the nilpotent pro-$p$
group~$\bfG_\Lambda(\mcO_\mfp)$. It is remarkable that many features
of the representation growth of groups of the form
$\bfG_\Lambda(\Gri_\mfp)$ only depend on the lattice $\Lambda$, and
not on the local ring $\Gri_\mfp$.

\begin{rem}
 We comment on connections between the above construction and the
 Mal'cev correspondence between $\T$-groups and nilpotent Lie
 rings. Starting from a $\T$-group $\Gamma$, there exists a $\Q$-Lie
 algebra $\mathcal{L}_\Gamma(\Q)$ and an injective mapping
 $\log:\Gamma\rarr\mathcal{L}_\Gamma(\Q)$, such that $\log(\Gamma)$
 spans $\mathcal{L}_\Gamma(\Q)$ over $\Q$. Whilst $\log(\Gamma)$ needs
 not, in general, be a Lie lattice inside $\mathcal{L}_\Gamma(\Q)$,
 there always exists a subgroup $H$ of $\Gamma$ of finite index with
 this property, satisfying $\log(H)'\subseteq c!\log(H)$, where $c$ is
 the nilpotency class of $\Gamma$. Setting~$\Lambda = \log(H)$, we
 recover $H$ as the group of $\Z$-rational points of~$\bfG_\Lambda$.
\end{rem}

Let now $\Lambda$ be again a nilpotent $\mcO$-Lie lattice of class
$c$, and suppose that $\La'\subseteq c!\La$. Denote by $\bfG_\Lambda$
the associated unipotent group scheme. For every finite extension $L$
of $K$, with ring of integers $\mcO_L$, we obtain a $\T$-group
$\bfG_\La(\mcO_L)$ and, for every nonzero prime ideal $\mfP\in\Spec
\mcO_L $, a pro-$p$ group~$\bfG_{\La}(\mcO_{L,\mfP})$.

\begin{thm}\cite{StasinskiVoll/14}\label{thm:SV_thm_A}
For every finite extension $L$ of $K$, with ring of integers $\mcO_L$,
\begin{equation}\label{equ:euler_T_groups}
  \zeta^{\wt{\textup{irr}}}_{\bfG_{\La}(\mcO_L)}(s) =
  \prod_{\mfP\in\Spec\mcO_L}\zeta^{\wt{\textup{irr}}}_{\bfG_\La(\mcO_{L,\mfP})}(s),
\end{equation} where, for each prime
ideal $\mfP\in\Spec\mcO_L$, the factor
$\zeta^{\wt{\textup{irr}}}_{\bfG_\La(\mcO_{L,\mfP})}(s)$ enumerates the continuous finite-dimensional irreducible
representations of
$\bfG_\La(\mcO_{L,\mfP})$ up to twisting by continuous one-dimensional
representations. Moreover, the following hold.

\begin{enumerate}
\item For each rational prime $p$,
\begin{equation*}
\zeta^{\wt{\textup{irr}}}_{\bfG_\La(\mcO_L),p}(s) =
\prod_{\mfP|p\mcO}\zeta^{\wt{\textup{irr}}}_{\bfG_\La(\mcO_{L,\mfP})}(s).
\end{equation*}

\item\label{rational} There exists a finite subset $S\subset \Spec
  \mcO$, an integer $t\in\N$, and a rational function
  $R(X_1,\dots,X_t,Y)\in\Q(X_1,\dots,X_t,Y)$ such that, for every
  prime ideal $\mfp\not\in S$, the following holds. There exist
  algebraic integers $\lambda_1,\dots,\lambda_t$, depending on $\mfp$,
  such that, for all finite extensions $\mfO$ of $\mfo=\mcO_\mfp$,
\begin{equation}\label{equ:rational.T}
\zeta^{\wt{\textup{irr}}}_{\bfG_\La(\mfO)}(s) =
R(\lambda_1^f,\dots,\lambda_t^f,q^{-fs}),
\end{equation}
where $q=|\mcO:\mfp|$ and $|\mcO_L:\mfP|=q^f$.

\item Setting $d=\dim_K(\Lambda'\otimes_\mcO K)$, the following
functional equation holds:
\begin{equation}\left.\zeta^{\wt{\textup{irr}}}_{\bfG_\La(\mfO)}(s)\right|_{\substack{q\rightarrow
      q^{-1}\\ \lambda_{i}\rightarrow\lambda_{i}^{-1}}} =
  q^{fd}\zeta^{\wt{\textup{irr}}}_{\bfG_\La(\mfO)}(s).
\end{equation}
\end{enumerate}
\end{thm}

As a corollary, we obtain that
$\zeta^{\wt{\textup{irr}}}_{\bfG_\La(\mfO)}(s)$ is rational in
$q^{-fs}$. In particular, the dimensions of the continuous
representations of the pro-$p$ group $\bfG_\La(\mfO)$ are all powers
of~$q^f$.

Example~\ref{exa:heisenberg.ztwirr} illustrates
Theorem~\ref{thm:SV_thm_A}. Indeed, the Heisenberg group scheme $\bfH$
is defined over $K=\Q$. We have $d=1$ and, in~\eqref{rational}, we may
take $S=\varnothing$, $t=1$, $R(X,Y)=\frac{1-Y}{1-XY}$ and
$\lambda_1=p$.

We say a few words about the proof of Theorem~\ref{thm:SV_thm_A},
referring to \cite{StasinskiVoll/14} for all details. The Euler
factorization~\eqref{equ:euler_T_groups} and the statement (1) follow
easily from strong approximation for unipotent groups. The key tool to
enumerate the representation zeta functions of pro-$p$ groups like
$\bfG_{\La}(\mfO)$ is the Kirillov orbit method. Wherever this method
is applicable, it parametrizes the irreducible representations of a
group in terms of the co-adjoint orbits in the Pontryagin dual of a
corresponding Lie algebra. In the case at hand, it reduces the problem
of enumerating twist-isoclasses of continuous finite-dimensional
irreducible representations of groups of the form $\bfG_{\La}(\mfO)$
to that of enumerating certain orbits in the duals of the derived
$\mfO$-Lie lattices~$\La(\mfO)'=(\La\otimes_\mcO\mfO)'$. By
translating the latter into the problem of evaluating $p$-adic
integrals, one reduces the problem further to the problem of
enumerating $p$-adic points on certain algebraic varieties, which only
depend on $\La$. In this way, one can show that there exist finitely
many smooth projective varieties $V_i$ defined over $\mcO$, and
rational functions $W_i(X,Y)\in\Q(X,Y)$, $i=1,\dots,N$, such that, if
$\mfp$ avoids a finite set $S\subset \Spec \Gri$,
$$\zeta^{\wt{\textup{irr}}}_{\bfG_\La(\mfO)}(s) = \sum_{i=1}^N
|\ol{V_i}(\mathbb{F}_{q^f})|W_i(q^f,q^{-fs}),$$ where $\ol{V_i}$
denotes reduction modulo~$\mfp$. By the Weil conjectures there exist,
for each $i\in\{1,\dots,N\}$, algebraic integers $\lambda_{ij}$,
$j=\{0,\dots,2\dim V_i\}$, such that
$$|\ol{V_i}(\mathbb{F}_{q^f})| = \sum_{j=0}^{2\dim V_i} (-1)^j
\lambda_{ij}^f$$ and $$\sum_{j=0}^{2\dim V_i} (-1)^j
\lambda_{ij}^{-f}= q^{f\dim W_i} \sum_{j=0}^{2\dim V_i} (-1)^j
\lambda_{ij}^f.$$ This remarkable symmetry is behind the functional
equations for the Hasse-Weil zeta functions of the
varieties~$\ol{V_i}$ and also functional equations such
as~\eqref{equ:funeq}. The rational functions $W_i$ come from the
enumeration of rational points of rational polyhedral cones.

\begin{qun}\label{qun:alpha.T}
  Given $\bfG_\La$ and $\mcO_L$ as above. Is the abscissa
  $\alpha^\text{irr}(\mathbf{G}_{\La}(\mcO_L))$ of
  $\zeta^{\wt{\textup{irr}}}_{\bfG_\La(\mcO_L)}(s)$ always a rational
  number? Is it independent of $L$?
\end{qun}

In general, the algebraic varieties $V_i$ are obtained from
resolutions of singularities of certain -- in general highly singular
-- varieties, and are difficult to compute explicitly. We give some of
the relatively few explicit examples of representation zeta functions
of $\T$-groups we have at the moment.

\begin{exm}
  Let $d\in\N_{>1}$ and $\mff_{d,2}$ the free nilpotent Lie ring on
  $d$ generators of nilpotency class $2$, of additive rank
  $d+\binom{d}{2}=\binom{d+1}{2}$. We write $\bfF_{d,2}$ for the
  unipotent group scheme $\bfG_{\mff_{d,2}}$ associated to this
  $\Z$-Lie lattice. For $d=2$ we obtain $\bfF_{2,2}=\bfH$, the
  Heisenberg group scheme. We also recover the free
  class-$2$-nilpotent group on $d$ generators as $\bfF_{d,2}(\Z)$. We
  write $d = 2\lfloor d/2 \rfloor + \varepsilon$ for $\varepsilon
  \in\{0,1\}$. The following
  generalizes~\eqref{equ:euler.ztwirr.Heisenberg.no.field}.

\begin{thm}\cite[Theorem~B]{StasinskiVoll/14}
  Let $\mcO$ be the ring of integers of a number field $K$. Then
$$\zeta^{\wt{\textup{irr}}}_{\bfF_{d,2}(\mcO)}(s) = \prod_{i=0}^{\lfloor d/2\rfloor} \frac{\zeta_K(s-2(\lfloor d/2 \rfloor + i + \varepsilon)+1)}{\zeta_K(s-2i)}.$$
\end{thm}
\end{exm}
E.~Avraham has computed the local factors of the representation zeta
function of the groups~$\bfF_{2,3}(\Gri[\frac{1}{6}])$; see
\cite{Avraham/13}. For further explicit examples of representation
zeta functions of $\T$-groups see \cite{Ezzat/12, Snocken/13}.

\subsection{Arithmetic lattices in semisimple groups}

Let $S$ be a finite set of places of a number field $K$, including all
archimedean ones, and let $\mcO_S$ denote the $S$-integers of~$K$. Let
further $\bfG$ be an affine group scheme over $\mcO_S$ whose generic
fibre is connected, simply-connected semi-simple algebraic group
defined over~$K$, together with a fixed embedding $\bfG
\hookrightarrow \GL_N$ for some~$N\in\N$. Let
$\Gamma=\bfG(\mcO_S)$. Then $\Gamma$ has polynomial representation
growth if and only if $\Gamma$ has the weak Congruence Subgroup
Property, i.e.\ the congruence kernel, that is the kernel of the
natural surjection
\begin{equation}\label{surjection}
\wh{\bfG(\mcO_S)} \rarr \bfG(\wh{\mcO_S})\cong
\prod_{\mfp\in(\Spec \mcO)\setminus S}\bfG(\mcO_\mfp),
\end{equation} 
is finite; cf.~\cite{LubotzkyMartin/04}. Here $\bfG(\wh{\mcO_S})$
denotes the congruence completion of~$\bfG(\mcO_S)$. For simplicity we
assume in the sequel that $\Gamma$ actually has the strong Congruence
Subgroup Property, i.e.\ that the congruence kernel is trivial, so
that the surjection~\eqref{surjection} is an isomorphism. A
prototypical example of such a group is the group $\SL_N(\Z)$
for~$N\geq 3$.

On the level of representation zeta functions, the triviality of the
congruence kernel is reflected by an Euler factorization, similar to
but different from those previously discussed, be it in the context of
subgroup and subring growth or of representation growth of
$\T$-groups. The Euler factorization features two types of factors:
the \emph{archimedean factors} are equal to
$\zeta^{\textup{irr}}_{\bfG(\C)}(s)$, the socalled \emph{Witten zeta
  function}, that is the Dirichlet generating series enumerating the
{rational} finite-dimensional irreducible complex representations of
the algebraic group $\bfG(\C)$. The \emph{non-archimedean factors}, on
the other hand, are the representation zeta functions
$\zeta^{\textup{irr}}_{\bfG(\mcO_\mfp)}(s)$, where $\mfp\not\in S$. These
Dirichlet generating series enumerate the {continuous}
finite-dimensional irreducible complex representations of the $p$-adic
analytic groups~$\bfG(\mcO_\mfp)$.

\begin{pro}\cite[Proposition~4.6]{LarsenLubotzky/08}
  The following Euler factorization holds:
\begin{equation}\label{equ:euler.arit}
\zeta^{\textup{irr}}_{\bfG(\mcO_S)}(s) = \zeta^{\textup{irr}}_{\bfG(\C)}(s)^{|K:\Q|}
\prod_{\mfp\in(\Spec \mcO)\setminus S}
\zeta^{\textup{irr}}_{\bfG(\mcO_\mfp)}(s).
\end{equation}
\end{pro}

It is a problem of central importance to compute the abscissa of
convergence $\alpha(\bfG(\mcO_S))$ of the representation zeta
function~$\zeta^{\textup{irr}}_{\bfG(\mcO_S)}(s)$. It is known that
$\alpha(\bfG(\mcO_S))$ is always a rational number; see
\cite[Theorem~1.2]{Avni/08} and compare Question~\ref{qun:alpha.T}.

The two types of factors of $\zeta^{\textup{irr}}_{\bfG(\mcO)}(s)$
in~\eqref{equ:euler.arit} turn out to have quite distinct flavours. We
discuss the archimedean local factors in
Section~\ref{subsubsec:witten}, the non-archimedean local factors in
\ref{subsubsec:p-adic}, and return to global zeta functions of
arithmetic groups in Section~\ref{subsubsec:global}.

\subsubsection{Witten zeta functions}\label{subsubsec:witten}
In this section let $\Gamma=\bfG(\C)$. For $n\in\N$ we denote by
$r_n(\Gamma)$ the number of $n$-dimensional rational, irreducible
complex representations of~$\Gamma$. Let $\Phi$ be the root system of
$\bfG$ of rank $r=\rk(\Phi)$, let $\Phi^+$ a choice of positive roots
of $\Phi$ and set $\rho = \sum_{\alpha\in\Phi^+}\alpha$. We write
$w_1,\dots,w_r$ for the fundamental weights. The rational irreducible
representations of $\Gamma$ are all of the form $W_\lambda$, where
$\lambda = \sum_{i=1}^r a_i w_i$ for~$a_i\in\N_0$. The Weyl dimension
formula asserts that
$$\dim W_\lambda = \prod_{\alpha\in\Phi^+}\frac{\langle
  \lambda+\rho,\alpha\rangle}{\langle \rho,\alpha \rangle}.$$ Note
that the numerator is a product of $\kappa:=|\Phi^+|$ affine linear
functions $f_1,\dots,f_\kappa$ in the integer coordinates of
$\lambda$, whilst the denominator $C=\prod_{\alpha\in\Phi^+}\langle
\rho,\alpha\rangle$ is a constant depending only on $\Phi$. Thus
\begin{equation}\label{equ:witten}
  \zeta^{\textup{irr}}_\Gamma(s) = \sum_\lambda (\dim W_\lambda)^{-s} = C^s \sum_{a\in\N_0^r}\prod_{i=1}^\kappa f_i(a)^{-s}.
\end{equation}

\begin{exm}
  Assume that $\bfG$ is of type $G_2$. Then $C=120$, $r=6$ and we may
  take
  \begin{alignat*}{3}
    f_1&=f_2=X_1+1,\quad &f_3 &= X_1+X_2+2,\quad& f_4&=X_1+2X_2+3,\\
    f_5 &= X_1 + 3X_2+4,\quad & f_6&=2X_1+3X_2+5 & &.
\end{alignat*}
\end{exm}

\begin{thm}\cite[Theorem~5.1]{LarsenLubotzky/08}\label{thm:abs.arch}
  The abscissa of convergence of $\zeta^{\textup{irr}}_{\bfG(\C)}(s)$
  is $r/\kappa$.
\end{thm}

Multivariable generalisations of zeta functions
like~\eqref{equ:witten} have been considered by Matsumoto
(\cite{Matsumoto/03}), among others. Functions of the form
$$\zeta(s_1,\dots,s_r;\bfG) = \sum_{a\in\N_0^r}\prod_{i=1}^r
f_i(a)^{-s_i},$$ where $s_1,\dots,s_r$ are complex variables, are, in
particular, known to have meromorphic continuation to the whole
complex plane; cf.~\cite[Theorem~3]{Matsumoto/03}.

Special values of Witten zeta functions are interpretable as volumes
of moduli spaces of certain vector bundles;
cf.~\cite[Section~7]{Zagier/92} and
\cite{Witten/91}. From~\eqref{equ:witten}, Zagier deduces
\begin{thm}\cite{Zagier/92} 
 If $s\in2\N$, then $\zeta^{\textup{irr}}_\Gamma(s)\in\Q\pi^{\kappa
   s}$.
\end{thm}

\subsubsection{Representation zeta functions of compact $p$-adic
  analytic groups}\label{subsubsec:p-adic} Let $\Gamma$ be a profinite
group. For $n\in\N$ we denote by $r_n(\Gamma)$ the number of
{continuous} finite-dimensional irreducible complex representations of
$\Gamma$. If $\Gamma$ is finitely generated, then $r_n(\Gamma)$ is
finite for all $n\in\N$ if and only if $\Gamma$ is FAb, i.e.\ has the
property that every open subgroup of $\Gamma$ has finite
abelianization.

\begin{thm}\cite[Theorem~1]{Jaikin/06}\label{thm:jaikin}
Let $p$ be an odd prime and $\Gamma$ a FAb compact $p$-adic analytic
group. Then there are natural numbers $n_1,\dots,n_k$ and rational
functions $W_1(Y),\dots,W_k(Y)\in\Q(Y)$ such that
\begin{equation}\label{equ:jaikin}
\zeta^{\textup{irr}}_\Gamma(s) = \sum_{i=1}^k n_i^{-s}W_i(p^{-s}).
\end{equation}
\end{thm}

\begin{exm}\label{exm:SL2}
Let $R$ be a compact discrete valuation ring whose (finite) residue
field $\Fq$ has odd characteristic. The representation zeta function
of the group $\SL_2(R)$ was computed in \cite[Section~7]{Jaikin/06}:
\begin{equation}\label{equ:SL2}
\zeta^{\textup{irr}}_{\SL_2(R)} (s) = \zeta^{\textup{irr}}_{\SL_2(\Fq)}(s)+ \left.\frac{
  4q\left(\frac{q^2-1}{2}\right)^{-s} + \frac{q^2-1}{2}(q^2-q)^{-s} +
  \frac{(q-1)^2}{2}(q^2+q)^{-s}}{1-q^{1-s}}\right.,
\end{equation} where
$$\zeta^{\textup{irr}}_{\SL_2(\Fq)}(s) = 1 + q^{-s} + \frac{q-3}{2}(q+1)^{-s} +
\frac{q-1}{2}(q-1)^{-s} + 2\left(\frac{q+1}{2}\right)^{-s} + 2\left(
\frac{q-1}{2}\right)^{-s}$$ is the representation zeta function of the
finite group of Lie type $\SL_2(\Fq)$.

If $R$ is a finite extension of $\Zp$, the ring of $p$-adic integers,
then \eqref{equ:SL2} illustrates~\eqref{equ:jaikin}. It is remarkable
that the same formula applies in the characteristic $p$ case, that is
if $R=\Fq\llbracket X \rrbracket$, the ring of formal power series
over $\Fq$.
\end{exm}
The proof of Theorem~\ref{thm:jaikin} utilizes the fact that a FAb
compact $p$-adic analytic group $\Gamma$ is virtually pro-$p$: it has
an open normal subgroup $N$ which one may assume to be uniformly
powerful. The Kirillov orbit method for uniformly powerful groups and
methods from model theory and the theory of definable $p$-adic
integration may be used to describe the distribution of the
representations of~$N$. Clifford theory is then applieded to extend
the analysis for $N$ to an analysis for $\Gamma$. The integers
$n_1,\dots,n_k$ are closely related to the dimensions of the
irreducible representations of the finite group $\Gamma/N$.

Computing zeta functions of FAb compact $p$-adic analytic groups --
such as the groups $\bfG(\mcO_\mfp)$ in \eqref{surjection} --
explicitly is in general very difficult. The situation is more
tractable for pro-$p$ groups. Theorem~\ref{thm:jaikin} states that if
$\Gamma$ is a FAb compact $p$-adic analytic pro-$p$ group, then
$\zeta^{\textup{irr}}_\Gamma(s)$ is rational in~$p^{-s}$. That this
generating function is a power series in $p^{-s}$ is obvious. Indeed,
the irreducible \emph{continuous} representations of a pro-$p$ group
$\Gamma$ all have $p$-power dimensions, as they factorize over finite
quotients of~$\Gamma$, which are all finite $p$-groups.

Representation zeta functions of pro-$p$ groups for which a version of
the Kirillov orbit method is available may be computed in terms of
$p$-adic integrals associated to polynomial mappings; see \cite[Part
  1]{AKOVI/13} for details. These integrals are of a much simpler type
than the general definable integrals used in the proof of
Theorem~\ref{thm:jaikin}.  In the following we discuss some cases
where this approach allows for an explicit computation of
representation zeta functions.

We concentrate on groups of the form $\bfG(\mfo)$, where $\mfo$ is a
finite extension of $\mcO_\mfp$ for some $\mfp\in(\Spec \mcO)\setminus
S$. Then $\mfo$ is a compact discrete valuation ring of characteristic
$0$, with maximal ideal $\mfm$, say, and finite residue field of
characteristic~$p$, where $\mfp | p\mcO$. For $m\in\N$ we consider the
\emph{$m$-th principal congruence subgroup} $\bfG^{m}(\mfo)$, that is
the kernel of the natural surjection
$$\bfG(\mfo) \rarr \bfG(\mfo/\mfm^m).$$ The groups $\bfG^{m}(\mfo)$
are FAb $p$-adic analytic pro-$p$ groups and, for sufficiently large
$m\in\N$, the Kirillor orbit method is applicable. This follows from
the fact that the groups $\bfG^m(\mfo)$ are saturable and potent for
$m \gg 0$; cf. \cite[Proposition~2.3]{AKOVI/13}
and~\cite{Gonzalez/09}. (In fact, if $\mfo$ is an unramified extension
of $\Zp$, then $m=1$ suffices.)  One would like to understand the
representation zeta functions
$\zeta^{\textup{irr}}_{\bfG^m(\mfo)}(s)$, and their variation with
\begin{itemize}
\item the prime ideal $\mfp\in(\Spec \mcO)\setminus S$, 
\item the ring extension $\mfo$, and
\item the congruence level~$m\in\N$.
\end{itemize}
The following result achieves much of this for the special linear
groups $\SL_3(\mfo)$ and the special unitary groups $\SU_3(\mfo)$,
assuming that $p\neq 3$. Here, the special unitary groups
$\SU_3(\mfo)$ are defined in terms of the nontrivial Galois
automorphism of the unramified quadratic extension of the field of
fractions of $\mfo$; see \cite[Section 6]{AKOVI/13} for details.

\begin{thm}\cite[Theorem~E]{AKOVI/13}\label{thm:SL3} Let $\mathfrak{o}$ be a
  compact discrete valuation ring of characteristic $0$ whose residue
  field has cardinality $q$ and characteristic different from~$3$. Let
  $\mathsf{G}(\mfo)$ be either $\SL_3(\mfo)$ or~$\SU_3(\mfo)$. Then,
  for all sufficiently large $m\in\N$,
\begin{equation}\label{equ:princ.A2}
    \zeta^{\textup{irr}}_{\mathsf{G}^m(\mfo)}(s) = q^{8m} \frac{1 + u(q) q^{-3-2s}
      + u(q^{-1}) q^{-2-3s} + q^{-5-5s}}{(1 - q^{1-2s})(1 - q^{2-3s})},
  \end{equation}
  where
  \begin{equation*} u(X) =
    \begin{cases}
      \phantom{-}X^3 + X^2 - X - 1 - X^{-1} & \text{ if }
      \mathsf{G}(\mfo)=\SL_3(\mfo),\\
      -X^3 + X^2 - X + 1 - X^{-1} & \text{ if }
      \mathsf{G(}\mfo)=\SU_3(\mfo).
    \end{cases}
  \end{equation*}
Furthermore, the following functional equation holds:
$$\left.\zeta^{\textup{irr}}_{\mathsf{G}^m(\mfo)}(s)\right|_{q \rarr
  q^{-1}} = q^{8(1-2m)}\zeta^{\textup{irr}}_{\mathsf{G}^m(\mfo)}(s).$$
\end{thm}

\begin{rem}\label{rem:SL3}
  We note that $\zeta^{\textup{irr}}_{\mathsf{G}^m(\mfo)}(s)$ is a
  rational function in $q^{-s}$ whose coefficients are given by
  polynomials in $q$, that $8$ is the dimension of the algebraic
  group~$\SL_3$, and that
  $\zeta^{\textup{irr}}_{\mathsf{G}^m(\mfo)}(s)/q^{8m}$ is independent
  of the congruence level $m$. Only a few signs in the numerators
  reflect the difference between special linear and unitary groups.
\end{rem}

In general, one can give formulae for the representation zeta
functions of groups of the form $\bfG^{m}(\mfo)$ --- valid for all
sufficiently large $m$ and virtually independent of $m$ --- which are
uniform both under variation of~$\mfp$ and~$\mfo$, and all but
independent of~$m$. More precisely, \cite[Theorem~A]{AKOVI/13} implies
the following result, which in turn generalizes Theorem~\ref{thm:SL3}.

\begin{thm}\cite[Theorem~A]{AKOVI/13}\label{thm:AKOV1_thm_A}
  There exist a finite subset $T \subset (\Spec \mcO) \setminus S$, an
  integer $t\in\N$, and a rational function
  $R(X_1,\dots,X_t,Y)\in\Q(X_1,\dots,X_t,Y)$ such that, for every
  prime ideal $\mfp\not\in S\cup T$, the following holds.

  There exist algebraic integers $\lambda_1,\dots,\lambda_t$,
  depending on~$\mfp$, such that, for all finite extensions $\mfO$ of
  $\mfo=\mcO_\mfp$, and all sufficiently large $m\in\N$,
\begin{equation}\label{equ:rational.arit}
\zeta^{\textup{irr}}_{\bfG^{m}(\mfO)}(s) =
q^{fdm}R(\lambda_1^f,\dots,\lambda_t^f,q^{-fs}),
\end{equation}
where $q=|\mcO:\mfp|$, $|\mfO:\mfP|=q^f$, and $d=\dim \bfG$.

Furthermore, the following functional equation holds:
\begin{equation}\label{equ:funeq.repgrow}
 \left.\zeta^{\textup{irr}}_{\bfG^{m}(\mfO)}(s)\right|_{\substack{q\rightarrow
     q^{-1}\\ \lambda_{i}\rightarrow\lambda_{i}^{-1}}} =
 q^{fd(1-2m)}\zeta^{\textup{irr}}_{\bfG^{m}(\mfO)}(s).
\end{equation}
\end{thm}

We note the close analogy between this result and
Theorem~\ref{thm:SV_thm_A}, which it precedes. Generalizing points
made in Remark~\ref{rem:SL3}, we further note that
Theorem~\ref{thm:AKOV1_thm_A} implies that
$\zeta^{\textup{irr}}_{\bfG^{m}(\mfO)}(s)$ is rational in~$q^{-fs}$
and $\zeta^{\textup{irr}}_{\bfG^{m}(\mfO)}(s)/q^{fdm}$ is independent
of~$m$. In general we do not expect that the coefficients of
$\zeta^{\textup{irr}}_{\bfG^{m}(\mfO)}(s)$ are given by polynomials in
$q^f$. In fact, as in Theorem~\ref{thm:SV_thm_A}, the algebraic
integers $\lambda_i$ arise from formulae for the numbers of rational
points of certain algebraic varieties over finite fields. One may ask,
however, whether these numbers are given by polynomials for
interesting classes of pro-$p$ groups arising from classical groups,
such as groups of the form $\SL_N^m(\mfo)$.

\begin{qun}\label{que:SLN}
  Let $N,m\in \N$ and $\mfo$ be a compact discrete valuation ring of
  characteristic $0$ whose residue field has cardinality $q$ and
  characteristic not dividing $N$. Does there exist a rational
  function $W_N(X,Y)\in\Q(X,Y)$ such that, for sufficiently large $m$,
$$\zeta^{\textup{irr}}_{\SL^m_N(\mfo)}(s) = q^{(N^2-1)m} W_N(q,q^{-s})?$$ The answer
  is ``yes'' in case $N=2$
  (cf.~\cite[Theorem~1.2]{AKOVIII/12}) and $N=3$
  (cf.\ Theorem~\ref{thm:SL3}).
\end{qun}

The striking similarity between the formulae for the representation
zeta functions of groups of the form $\SL_3^m(\mfo)$ and $\SU_3(\mfo)$
is reminiscent of Ennola duality for the characters of the finite
groups $\GL_n(\Fq)$ and $\GU_n(\Fq)$; cf.~\cite{Kawanaka/85}. I am not
aware of such a duality in the realm of compact $p$-adic analytic
groups, but read \eqref{equ:princ.A2} as a strong indication for a
connection like this.

Computing the representation zeta functions of the ``full'' $p$-adic
analytic groups $\bfG(\mcO_\mfp)$ is significantly harder than those
of their principal congruence subgroups. In principle, Clifford theory
allows one to describe the representations of the former groups in
terms of the representations of their open normal subgroups
$\bfG^m(\mcO_\mfp)$. However, how to tie in explicit Clifford theory
with the theory that leads to results like
Theorem~\ref{thm:AKOV1_thm_A} in a way that is uniform in $\mfp$ and
$\mfo$ is not clear in general.

The paper \cite{AKOVII/14} contains formulae for the representation
zeta functions of special linear groups of the form $\SL_3(\mfo)$ and
special unitary groups of the form $\SU_3(\mfo)$, where $\mfo$ is an
unramified extension of $\Zp$ and~$p\neq 3$.  The resulting formulae
of the form~\eqref{equ:jaikin} are significantly more complicated than
the formulae \eqref{equ:princ.A2} for the principal congruence
subgroups, and are omitted here. We just record the fact that
$$(1-q^{1-2s})(1-q^{2-3s})$$ 
is a common denominator for the rational functions involved, just as
in~\eqref{equ:princ.A2}.

It is of great interest if these formulae also apply in characteristic
$p$, i.e.\ for groups like $\SL_3(\Fq\llbracket X \rrbracket)$. In
contrast to the hands-on computations in \cite{Jaikin/06}, the
computations in \cite{AKOVII/14} do rely on the Kirillov orbit method
for uniformly powerful subgroups of the relevant $p$-adic analytic
groups, which is only available in characteristic~$0$.
 
In \cite{AKOVII/14} we also compute the representation zeta functions
of finite quotients of groups of the form
\begin{gather*}
\SL_3(\mfo), \SU_3(\mfo), \GL_3(\mfo),\GU_3(\mfo), \\ \SL_3^m(\mfo),
\SU_3^m(\mfo), \GL_3^m(\mfo),\GU_3^m(\mfo)
\end{gather*}
by principal congruence subgroups, subject to some restrictions on the
residue field characteristic~$p$. Some further examples of
representation zeta functions of $p$-adic analytic groups are
contained in \cite{AKOVIII/12}. Recent results of Aizenbud and Avni
bound the abscissa of convergence of zeta functions of groups of the
form $\SL_N(\lri)$, where $\lri$ is a compact discrete valuation ring
of characteristic~$0$; cf.~\cite[Theorem~A]{AizenbudAvni/13}.

We close this section by mentioning a vanishing theorem for
representation zeta functions.

\begin{thm}\cite{GonzalezJaikinKlopsch/14}
  Let $p$ be an odd prime and $\Gamma$ an infinite FAb compact
  $p$-adic analytic group. Then $\zeta^{\textup{irr}}_\Gamma(-2)=0$.
\end{thm}

The proof of this result uses the fact that, while the series
$\zeta^{\textup{irr}}_\Gamma(s)$ does not converge in the usual topology for
$s\in\R_{<0}$, the expressions $\zeta^{\textup{irr}}_\Gamma(e)$ do converge in the
$p$-adic topology for all negative integers~$e$.

\subsubsection{Representation zeta functions of arithmetic lattices}
\label{subsubsec:global}
We now return to the global representation zeta function of
${\bfG(\mcO_S)}$.

For the purpose of analyzing $\zeta^{\textup{irr}}_{\bfG(\mcO_S)}(s)$
via the Euler factorization \eqref{equ:euler.arit}, uniform formulae
for zeta functions of the form
$\zeta^{\textup{irr}}_{\bfG^m(\mcO_\mfp)}(s)$ --- as provided, e.g.,
by Theorem~\ref{thm:AKOV1_thm_A} --- are of limited value. Indeed,
whilst the index of $\bfG^m(\mcO_\mfp)$ in $\bfG(\mcO_\mfp)$ is finite
for each $\mfp$ and all $m$, the representation zeta function of every
finite index subgroup of $\bfG(\mcO)$ will share all but finitely many
of its non-archimedean factors with those of
$\zeta^{\textup{irr}}_{\bfG(\mcO_S)}(s)$.

Essentially only for groups of type $\mathsf{A}_2$ do we know how to
use Clifford effectively to deduce explicit uniform formulae for the
representation zeta functions
$\zeta^{\textup{irr}}_{\bfG(\mcO_\mfp)}(s)$;
cf.~\cite{AKOVII/14}. This allows for precise asymptotic results about
the representation growth of arithmetic groups of type $\mathsf{A}_2$.

\begin{thm}\cite{AKOVII/14}\label{thm:asy.A2}
  Let $\bfG$ be a connected, simply-connected absolutely simple
  algebraic group defined over $K$ of type $\mathsf{A}_2$, and assume
  that $\Gamma = \bfG(\Gri_S)$ has the strong Congruence Subgroup
  Property. Then $\alpha(\Gamma)=1$. Moreover,
  $\zeta^{\textup{irr}}_\Gamma(s)$ admits meromorphic continuation to
  $\{s\in\C \mid \Re(s) > 5/6\}$. The continued function is analytic
  on this half-plane, except for a double pole at $s=1$. Consequently,
  there exists a constant $c(\Gamma)\in\R_{>0}$, such that
$$\sum_{i=1}^n r_i(\Gamma) \sim c(\Gamma) \cdot  n \log n.$$
\end{thm}

We comment briefly on the proof of Theorem~\ref{thm:asy.A2}. Let
$\Gamma$ be as in the theorem. It is a key fact that all but finitely
many of the Euler factors of $\zeta^{\textup{irr}}_\Gamma(s)$ are of
the form $\SL_3(\mcO_\mfp)$ or $\SU_3(\mcO_\mfp)$, where $\mfp$ is a
prime ideal of~$\Gri$. To see that $\alpha(\Gamma)=1$, it suffices to
prove that the abscissa of convergence of the product over these
factors is equal to $1$. Indeed, \cite[Theorem~B]{AKOVI/13} implies
that the abscissa of convergence of the Euler
factorization~\eqref{equ:euler.arit} remains unchanged by removing
finitely many non-archimedean factors. The archimedean factors'
abscissa of convergence is~$2/3$; cf.~Theorem~\ref{thm:abs.arch}. To
compute the abscissa of convergence of the Euler factorization of the
factors of the form $\SL_3(\Gri_\mfp)$ or $\SU_3(\Gri_\mfp)$, one may
either inspect the explicit formulae given in~\cite{AKOVII/14}, or
argue with ``approximative Clifford theory'' as in \cite{AKOVI/13}.

The existence of meromorphic continuation is evident from inspection
of the explicit formulae for
$\zeta^{\textup{irr}}_{\SL_3(\mcO_\mfp)}(s)$ and
$\zeta^{\textup{irr}}_{\SU_3(\mcO_\mfp)}(s)$. The key here is that the
relevant Euler factorization can be approximated by the following
product of translates of (partial) Dedekind zeta functions:
\begin{equation}\label{equ:dedekind}
  \zeta_{K,S}(2s-1)\zeta_{K,S}(3s-2) = \prod_{\mfp\in(\Spec \mcO)\setminus S} \frac{1}{(1-| \mcO : \mfp|^{1-2s})(1-| \mcO : \mfp|^{2-3s})}.
\end{equation}
Roughly speaking, dividing $\zeta^{\textup{irr}}_{\SL_3(\mcO_\mfp)}(s)$ or
$\zeta^{\textup{irr}}_{\SU_3(\mcO_\mfp)}(s)$ by the appropriate local factor
of~\eqref{equ:dedekind} clears their common denominator
$(1-q^{1-2s})(1-q^{2-3s})$, and the Euler factorization of the remaining
numerators converges strictly better than the original Euler factorization.

Theorem~\ref{thm:asy.A2} states, in particular, that the abscissa of
convergence of the representation zeta function of an arithmetic group
of type $\mathsf{A}_2$ is always equal to $1$: the degree of representation
growth of very different groups --- such as, for example,
$\SL_3(\mcO)$ and $\SU_3(\mcO)$, for various number rings $\mcO$ ---
only depends on the root system of the underlying algebraic
group. This remarkable fact is vastly generalized by the following
result.

\begin{thm}\cite[Theorem~1.1]{AKOVIV/14}\label{thm:bc}
  Let $\Phi$ be an irreducible root system. Then there exists a
  constant $\alpha_\Phi\in\Q$ such that, for every arithmetic group
  $\bfG(\mcO_S)$, where $\mcO_S$ is the ring of $S$-integers of a
  number field $K$ with respect to a finite set of places $S$ and
  $\bfG$ is a connected abolutely almost simple algebraic group over
  $K$ with absolute root system $\Phi$, the following holds: if
  $\bfG(\mcO_S)$ has the CSP, then $\alpha(\bfG(\mcO_S)) =
  \alpha_\Phi$.
\end{thm}

Theorem~\ref{thm:bc} reduces a conjecture of Larsen and Lubotzky on
the invariance of representation growth of lattices in higher rank
semisimple locally compact groups to a conjecture of Serre on the CSP;
see \cite[Theorem~1.3]{AKOVIV/14}. A key idea of its proof is to
approximate the local factors of the representation zeta function
$\zeta^{\textup{irr}}_{\bfG(\mcO_S)}(s)$ uniformly by certain
definable integrals, in a way that leaves the abscissa of convergence
unchanged. The proof uses deep, nonconstructive techniques from model
theory, which hold little promise to yield an explicit description of
the function $\Phi\mapsto \alpha_\Phi$. So far, the only explictly
known values of this function are $\alpha_{A_1}=2$ and
$\alpha_{\mathsf{A}_2}=1$. 

\begin{qun}
  What is the value of $\alpha_\Phi$ in Theorem~\ref{thm:bc}, for
  various root systems $\Phi$?
\end{qun}

\subsection{Iterated wreath products and branch groups}

Let $Q$ be a finite group, acting on a finite set $X$ of cardinality
$|X|=d\geq 2$. We define iterated permutational wreath products as
follows. Set $W(Q,0):=\{1\}$ and, for $k\in\N$, set $W(Q,k+1) =
W(Q,k)\wr_X Q$. Passing to the inverse limit yields the profinite
group $W(Q) := \varprojlim_kW(Q,k)$. Recall that, for a profinite
group $G$, we denote by $r_n(G)$ the number of continuous
$n$-dimensional irreducible complex representations of $G$ up to
isomorphism and that $G$ is called rigid if $r_n(G)<\infty$ for all
$n\in\N$.

\begin{thm}\cite{BartholdidelaHarpe/10}
$W(Q)$ is rigid if and only if the group $Q$ is perfect,
  i.e.\ $G=[G,G]$. In this case, the following hold.
\begin{enumerate}
 \item The abscissa of convergence
   $\alpha:=\alpha(\zeta^{\textup{irr}}_{W(Q)}(s))$ is positive and
   finite, i.e.\ $\alpha\in\R_{>0}$.
 \item Locally around $\alpha$, the function
   $\zeta^{\textup{irr}}_{W(Q)}(s)$ allows for a Puiseux expansion of
   the form $$\sum_{n=0}^\infty c_n(s-\alpha)^{n/e}$$ for suitable
   $c_n\in\C$, $n\in\N$, and $e\in\{2,3,\dots,d\}$.
 \item Let $p_1,\dots,p_\ell$ denote the primes dividing $|Q|$. There
   exists a nontrivial polynomial
   $\Psi\in\Q[X_1,\dots,X_d,Y_1,\dots,Y_\ell]$ such that
  \begin{equation}\label{equ:recursion}
    \Psi(\zeta^{\textup{irr}}_{W(Q)}(s),\zeta^{\textup{irr}}_{W(Q)}(2s),\dots,\zeta^{\textup{irr}}_{W(Q)}(ds),p_1^{-s},\dots,p_\ell^{-s})=0.
 \end{equation}
 \end{enumerate}
 \end{thm}

 For examples illustrating in particular the functional
 equations~\eqref{equ:recursion}, see~\cite{BartholdidelaHarpe/10}.
 For generalizations of these results to self-similar profinite
 branched groups, see~\cite{Bartholdi/13}.

\begin{acknowledgements}
  Over several years, my work has been generously supported by
  numerous funding bodies, including the EPSRC, the DFG, and the
  Nuffield Foundation. I am also grateful to the organizers of Groups
  St Andrews 2013 in St~Andrews.
\end{acknowledgements}


\def\cprime{$'$}
\providecommand{\bysame}{\leavevmode\hbox to3em{\hrulefill}\thinspace}
\providecommand{\MR}{\relax\ifhmode\unskip\space\fi MR }
\providecommand{\MRhref}[2]{%
  \href{http://www.ams.org/mathscinet-getitem?mr=#1}{#2}
}
\providecommand{\href}[2]{#2}

\end{document}